\documentclass[12pt]{article}
\usepackage{amsmath}
\usepackage{amssymb}
\usepackage{latexsym}
\usepackage{amsthm}
\usepackage{tikz}
\usepackage{tabularx}
\usepackage{booktabs}
\usepackage{mathrsfs}
\usepackage{enumitem}
\usetikzlibrary{calc}
\usepackage{xcolor} 
\usepackage{ytableau}
\usepackage{graphicx}
\usepackage{algorithm,algorithmic}
\usepackage[colorlinks,
linkcolor=red,
anchorcolor=blue,
citecolor=green]{hyperref}

\newtheorem{thm}{Theorem}[section]
\newtheorem{lem}[thm]{Lemma}

\newtheorem{cor}[thm]{Corollary}
\newtheorem{conj}[thm]{Conjecture}
\newtheorem{rem}[thm]{Remark}

\begin{document}
	\begin{center}
		{\large \bf The Eulerian transformation preserves log-convexity}
	\end{center} 
\begin{center}
	Huihua Gao$^1$, Xin-Bei Liu$^2$, Weilin Wang$^3$ \\[6pt]

	$^{1,2}$Center for Combinatorics, LPMC, Nankai University,\\
	Tianjin 300071, P. R. China\\
	$^3$School of Mathematical Sciences, Qufu Normal University,\\
	Qufu 273165, P. R. China\\[6pt]

	Email: $^1${\tt huihuagao@mail.nankai.edu.cn},
	$^2${\tt lxb@mail.nankai.edu.cn},\\
	$^3${\tt wangweilin@qfnu.edu.cn}
\end{center}
	\noindent\textbf{Abstract.}
In their 2007 paper on the log-convexity of combinatorial sequences,     
Liu and Wang conjectured that the Eulerian transformation preserves log-convexity. We completely solve their conjecture by showing that this transformation preserves strong $q$-log-convexity. More generally, we provide a suﬀicient
condition for a linear transformation, defined by a triangular 
array satisfying a certain three-term recurrence relation, to preserve strong $q$-log-convexity. 
As applications, we show that both the $r$-colored Eulerian transformation and
the $m$-Stirling Eulerian transformation preserve strong $q$-log-convexity. In particular, we confirm a conjecture due to Liu and Zhu, which states that the Eulerian transformation of type $B$ preserves log-convexity. 


\noindent \emph{AMS Classification 2020:}  05A05; 05A20; 26D15

\section{Introduction}

The main objective of this paper is to prove a challenging conjecture due to Liu and Wang \cite[Conjecture~5.4]{Liu-Wang-2007}, which states that the Eulerian transformation preserves log-convexity. We provide two proofs of Liu and Wang's conjecture, one of which relies on a careful use of the Cauchy–Schwarz inequality, and the other establishes the strong $q$-log-convexity preserving property of this transformation. Each proof provides a unified approach to preserving the log-convexity or strong 
$q$-log-convexity of many linear transformations.   

Let us first give an overview of related definitions and background. Recall that a sequence $(x_n)_{n\ge0}$ of nonnegative real numbers is called \emph{log-convex} (resp.\ \emph{log-concave}) if  
$x_n^2\le x_{n-1}x_{n+1}$ (resp.\ $x_n^2\ge x_{n-1}x_{n+1}$) for every $n\ge1$. We say that it is \emph{strongly log-convex} (resp.\ \emph{strongly log-concave}) if 
$x_nx_{n+m+1}\ge x_{n+1}x_{n+m}$
(resp.\ $x_nx_{n+m+1}\le x_{n+1}x_{n+m}$)
for all $n,m\ge0$.  Clearly,  
the log-convexity (resp.\ log-concavity) of a nonnegative sequence with no internal zeros is equivalent to its strong log-convexity (resp.\ log-concavity). 
The log-convex and log-concave sequences have attracted considerable attention in
combinatorics, algebra, geometry, analysis, probability and statistics; see \cite{Brascamp-Lieb-1976,Brenti-1995,Davenport-Polya-1949,Liu-Wang-2007,Saumard-Wellner-2014,Stanley-1989} for instance.
This paper is mainly concerned with log-convexity. 
Generally, it is not easy to prove the log-convexity of certain sequences.
Sometimes, we may fulfill this aim by using certain linear transformations associated with triangular arrays that preserve log-convexity. 
Given a triangular array
$T=\{T_{n,k}\}_{n,k\geq0}$ of nonnegative real numbers and $T_{n,k}=0$ unless $0\le k\le n$, we say that the linear transformation 
$z_n=\sum_{k=0}^{n}T_{n,k}x_k$ preserves log-convexity if
$(z_n)_{n\geq0}$ is log-convex whenever $(x_n)_{n\geq0}$ is a
nonnegative log-convex sequence.  
The preservation of log-convexity by several classical combinatorial transformations has been established. The result for the binomial transformation follows from a theorem of Davenport and P\'olya~\cite{Davenport-Polya-1949} on binomial convolution. Liu and Wang~\cite{Liu-Wang-2007} proved the corresponding results for the Stirling transformations of both kinds, while Liu~\cite{Liu-2011} established the same property for the Lah and Dowling transformations. The Bessel transformation was shown to preserve log-convexity by Chen, Wang and Yang~\cite{Chen-Wang-Yang-2011}. The preservation of log-convexity by the Narayana transformations of two kinds was conjectured by Liu and Wang~\cite{Liu-Wang-2007}. These conjectures were subsequently proved by Chen, Wang and Yang~\cite{Chen-Wang-Yang-2010} and Chen, Tang, Wang and Yang~\cite{Chen-Tang-Wang-Yang-2010}, respectively.

This paper is motivated by the log-convexity-preserving conjecture of the Eulerian transformation, which was proposed by Liu and Wang \cite{Liu-Wang-2007}. For $n\ge1$, let $[n]=\{1,\ldots,n\}$ and let
$\mathfrak{S}_n$ denote the set of permutations of $[n]$. Recall that the Eulerian number $A(n,k)$ counts permutations
$w=w_1\cdots w_n\in\mathfrak{S}_n$ with exactly $k-1$
descents, where a descent is an index $1\le i<n$
such that $w_i>w_{i+1}$. It is known that the Eulerian numbers satisfy
\begin{align*}
A(n,k)=kA(n-1,k)+(n-k+1)A(n-1,k-1),
\quad 0\le k\le n,
\end{align*}
with $A(0,0)=1$ and $A(n,k)=0$ for $k<0$ or $k>n$. 
We refer the transformation $z_n=\sum_{k=0}^{n}A(n,k)x_k$ as the 
Eulerian transformation, which differs from the notion of the same name considered by Athanasiadis~\cite{Athanasiadis-2025}, Br\"anden and Jochemko \cite{Branden-Jochemko-2022}, and Brenti \cite{Brenti-1989}. Liu and Wang's conjecture is stated as follows. 

\begin{conj}[{\cite[Conjecture~5.4]{Liu-Wang-2007}}]\label{main_conj}
If $(x_n)_{n\geq0}$ is a log-convex sequence, so is $(z_n)_{n\geq0}$, where $z_n=\sum_{k=0}^{n}A(n,k)x_k.$
\end{conj}

Some efforts have been made to prove the above conjecture. In particular, Wang and Zhu~\cite{Wang-Zhu-2016} proved that the Eulerian transformation preserves the Stieltjes moment property, a stronger property than log-convexity. 
Recall that a matrix of nonnegative real numbers is called \emph{totally positive of order $r$}
($\mathrm{TP}_r$ for short) if all its minors of order at most $r$ are nonnegative.
It is called \emph{totally positive} (TP for short) if it is $\mathrm{TP}_r$ for every $r\ge1$.
For a nonnegative sequence $x=(x_n)_{n\ge0}$, let
$\mathcal{H}(x)=[x_{i+j}]_{i,j\ge0}$ be its Hankel matrix.
Then $x$ is log-convex if and only if $\mathcal{H}(x)$ is
$\mathrm{TP}_2$. We say that $x$ is a \emph{Stieltjes moment sequence}
(or an \emph{SM sequence}) if $\mathcal{H}(x)$ is TP. The preservation of the SM property under linear transformations and convolutions was studied by Wang and Zhu~\cite{Wang-Zhu-2016}, who provided a sufficient condition with applications to the Stirling transformation of the second kind, the Narayana transformations of two kinds, and the Eulerian transformation. 
However, Wang and Zhu's result does not imply Conjecture \ref{main_conj}, since a nonnegative log-convex sequence need not be an SM sequence. 

Our first approach to Conjecture \ref{main_conj} involves some triangular array $T = \{T_{n,k}\}_{n,k\ge0}$ of nonnegative real numbers satisfying the recurrence 
\begin{align}\label{T_recurrence}
T_{n+1,k}
=\alpha_{n,k}T_{n,k}+\beta_{n,k-1}T_{n,k-1},
\end{align}
where $\alpha_{n,k}, \beta_{n,k-1}$ are real numbers, $T_{0,0}=1$ and $T_{n,k}=0$ unless $0\le k\le n$.
We have the following result. 

\begin{thm}\label{log-convex}
Let $T=\{T_{n,k}\}_{n,k\ge0}$ be a triangular array
of nonnegative real numbers satisfying
\eqref{T_recurrence}.
Assume that, for $n\ge0$ and $0\le k\le n$,
\begin{align*}
\alpha_{n,k}&\ge0,
&\alpha_{n+1,k}&\ge\alpha_{n,k},
&\alpha_{n+1,k+1}&\ge\alpha_{n,k},\\
\beta_{n,k}&\ge0,
&\beta_{n+1,k}&\ge\beta_{n,k},
&\beta_{n+1,k+1}&\ge\beta_{n,k}.
\end{align*}
Then the linear transformation
$z_n=\sum_{k=0}^nT_{n,k}x_k$ preserves log-convexity.
\end{thm}

We are able to prove Conjecture \ref{main_conj} 
as a special case of Theorem \ref{log-convex} by taking $\alpha_{n,k}=k$ and $\beta_{n,k}=n+1-k$ in \eqref{T_recurrence}. It is worth mentioning that Liu and Pei~\cite{Liu-Pei} independently confirmed this conjecture by proving the following stronger statement: if both $(x_n)_{n\geq0}$ and $(y_n)_{n\geq0}$ are log-convex, so is 
the sequence $(z_n)_{n\geq0}$, where $z_n=\sum_{k=0}^{n}A(n,k)x_k y_{n-k}$.

The second approach to Conjecture \ref{main_conj} of this paper is to show that the Eulerian transformation preserves strong $q$-log-convexity, a polynomial analogue of strong log-convexity. 
Let $\mathbb R_{\ge0}[q]$ denote the set of polynomials
in $q$ with nonnegative real coefficients.
For real polynomials $f(q)$ and $g(q)$, write
$f(q)\ge_q g(q)$ if
$f(q)-g(q)\in\mathbb R_{\ge0}[q]$.
A sequence $(X_n(q))_{n\ge0}$ in $\mathbb R_{\ge0}[q]$
is called \emph{$q$-log-convex}~\cite{Liu-Wang-2007} if
\begin{align*}
\mathcal{T}_{n}(X;q):=X_{n-1}(q)X_{n+1}(q)-X_n(q)^2\ge_q0,
\qquad n\ge1,
\end{align*}
and \emph{strongly $q$-log-convex}~\cite{Chen-Wang-Yang-2011} if
\begin{align*}
\mathcal{T}_{n,m}(X;q):=X_n(q)X_{n+m+1}(q)-X_{n+1}(q)X_{n+m}(q)\ge_q0,
\qquad n,m\ge0.
\end{align*}
The corresponding notion of $q$-log-concavity
(resp.\ strong $q$-log-concavity), introduced by
Stanley~\cite{Stanley-1989}
(resp.\ Sagan~\cite{Sagan-1992}),
is obtained by reversing the respective inequality.
The strong $q$-log-convexity implies $q$-log-convexity,
although the converse need not hold.
Similarly, a matrix with entries in $\mathbb R[q]$
is called \emph{$q$-totally positive}
($q$-TP for short) if all its minors have nonnegative
coefficients.
A sequence $(X_n(q))_{n\ge0}$ in $\mathbb R_{\ge0}[q]$
is called a \emph{$q$-Stieltjes moment sequence}
(or a \emph{$q$-SM sequence}) if its Hankel matrix
$[X_{i+j}(q)]_{i,j\ge0}$ is $q$-TP;
see Wang and Zhu~\cite{Wang-Zhu-2016}.
The $q$-SM property implies strong $q$-log-convexity,
which in turn implies $q$-log-convexity.

As mentioned earlier, the Eulerian transformation preserves the SM property. 
In this paper we shall prove that 
the Eulerian transformation preserves strong $q$-log-convexity. 
The transformation associated with the triangular array
$T=\{T_{n,k}\}_{n,k\geq0}$ is said to preserve strong
$q$-log-convexity if the strong $q$-log-convexity of
$(X_n(q))_{n\geq0}$ implies that of $(Z_n(q))_{n\geq0}$, where $Z_n(q)=\sum_{k=0}^{n}T_{n,k}X_k(q)$.
Several sufficient conditions for linear transformations to preserve strong
 $q$-log-convexity have been obtained.
Zhu and Sun~\cite{Zhu-Sun-2015} gave a sufficient condition
for preserving strong $q$-log-convexity and applied it to the binomial,
Morgan--Voyce and Narayana transformations of two kinds.
Liu and Li~\cite{Liu-Li-2016} obtained a criterion based on triangular
recurrences, with applications to the Stirling and Whitney
transformations. Further criteria and applications to the
Jacobi--Stirling, central factorial and Ramanujan transformations
were obtained by Zhu~\cite{Zhu-2018}. But neither of these known
 criteria applies to the Eulerian transformation.

As in the first approach to Conjecture \ref{main_conj}, we study the strong $q$-log-convexity preserving property of the Eulerian transformation 
in a more general framework. We specialize the coefficients in
\eqref{T_recurrence} to
$\alpha_{n,k}=a_n+c_nk$ and $\beta_{n,k}=b_n-c_nk$,
where $(a_n)_{n\geq0}$, $(b_n)_{n\geq0}$, and
$(c_n)_{n\geq0}$ are real sequences. Thus
\eqref{T_recurrence} becomes
\begin{align}\label{T_strong}
T_{n+1,k}
=(a_n+c_nk)T_{n,k}
+\bigl(b_n-c_n(k-1)\bigr)T_{n,k-1},
\qquad 0\leq k\leq n+1.
\end{align}
We have the following result. 

\begin{thm}\label{thm:strong-q-log-convexity-preserving}
Let $T=\{T_{n,k}\}_{n,k\geq0}$ be a triangular array of
nonnegative real numbers satisfying \eqref{T_strong}, with
$T_{0,0}=1$. Suppose that $a_0,\,b_0\geq0$,  and
\begin{align}
a_{n+1}\geq a_n,\quad
c_{n+1}\geq c_n\geq0,\quad
b_{n+1}-b_n\geq (n+1)c_{n+1}-nc_n,
\quad n\geq0. \label{eq-thm3-condition}
\end{align}
Then the linear transformation
$Z_n(q)=\sum_{k=0}^{n}T_{n,k}X_k(q)$
preserves strong $q$-log-convexity.
\end{thm}

Note that if $T=\{T_{n,k}\}_{n,k\geq0}$ satisfies the conditions of the above theorem, then it also satisfies the hypotheses of
Theorem~\ref{log-convex}. 
Taking $a_n=0$, $b_n=n+1$, and $c_n=1$ for all $n\geq0$
in Theorem~\ref{thm:strong-q-log-convexity-preserving}, we
obtain the following corollary, which strengthens Liu and Wang's conjecture.

\begin{cor}\label{cor:eulerian-strong-q-log-convexity-preserving}
The Eulerian transformation
$Z_n(q)=\sum_{k=0}^{n}A(n,k)X_k(q)$
preserves strong $q$-log-convexity.
\end{cor}

The paper is organized as follows. In 
Section~\ref{Log-Convexity Preservation} we give a proof of
Theorem~\ref{log-convex} based on the Cauchy--Schwarz inequality.
In Section~\ref{Strong $q$-Log-Convexity Preservation} we prove
Theorem~\ref{thm:strong-q-log-convexity-preserving}.
In Section~\ref{sec:applications}, we present several applications of
Theorems~\ref{log-convex} and
\ref{thm:strong-q-log-convexity-preserving}. In particular, we show
that both the $r$-colored Eulerian transformation and the
$m$-Stirling Eulerian transformation preserve strong
$q$-log-convexity.

\section{Log-convexity preservation}\label{Log-Convexity Preservation}

The main objective of this section is to prove Theorem~\ref{log-convex}. 
To this end, we first recall the following Cauchy--Schwarz inequality, which plays an important role in our proof of Theorem~\ref{log-convex}.

\begin{lem}[{\cite[Chapter~1]{Steele-2004}}]\label{cauchy-schwarz}
Let $s\geq0$ be a nonnegative integer.
For any real numbers $\xi_0,\ldots,\xi_s$ and
$\eta_0,\ldots,\eta_s$, we have
\begin{align*}
\left(\sum_{k=0}^s \xi_k^2\right)
\left(\sum_{k=0}^s \eta_k^2\right)
-\left(\sum_{k=0}^s \xi_k\eta_k\right)^2
\geq0.
\end{align*}
\end{lem}
We now proceed to prove Theorem~\ref{log-convex}.
\begin{proof}[Proof of Theorem~\ref{log-convex}]
We first assume that $x_k>0$ for every $k\geq0$. Set $\rho_k={x_{k+1}}/{x_k}$ for $k\geq0$.
Since the sequence $(x_k)_{k \ge 0}$ is log-convex, its quotient sequence $(\rho_k)_{k\ge0}$ is weakly increasing.
For $0\leq k\leq n$, set
\begin{align*}
w_{n,k}=T_{n,k}x_k,\quad 
g_{n,k}=\alpha_{n,k}+\beta_{n,k}\rho_k,\quad  \ell_{n,k}=\alpha_{n,k}g_{n+1,k}
+\beta_{n,k}\rho_k g_{n+1,k+1}.   
\end{align*}
Using the recurrence \eqref{T_recurrence}, we obtain
\begin{align}\label{eq:z-n-plus-1}
z_{n+1}
&=\sum_{k=0}^{n+1}T_{n+1,k}x_k=\sum_{k=0}^{n+1}\alpha_{n,k}T_{n,k}x_k+\sum_{k=0}^{n+1}\beta_{n,k-1}T_{n,k-1}x_k \nonumber \\
&=\sum_{k=0}^{n}\alpha_{n,k}T_{n,k}x_k+\sum_{k=0}^{n}\beta_{n,k}T_{n,k}x_{k+1}\nonumber\\
&=\sum_{k=0}^{n}T_{n,k}x_k\bigl(\alpha_{n,k}+\beta_{n,k}\rho_{k}\bigr)\nonumber\\
&=\sum_{k=0}^{n}w_{n,k}g_{n,k}.
\end{align}
In turn, one can verify that 
\begin{align}\label{eq:z-n-plus-2}
z_{n+2}
&=\sum_{k=0}^{n+1}T_{n+1,k}x_kg_{n+1,k}=\sum_{k=0}^{n+1}(\alpha_{n,k}T_{n,k}+\beta_{n,k-1}T_{n,k-1})x_kg_{n+1,k}\nonumber\\
&=\sum_{k=0}^{n}w_{n,k}
\bigl(\alpha_{n,k}g_{n+1,k}
+\beta_{n,k}\rho_k g_{n+1,k+1}\bigr)\nonumber\\
&=\sum_{k=0}^{n}w_{n,k}\ell_{n,k}
\end{align}
Combining \eqref{eq:z-n-plus-1} and
\eqref{eq:z-n-plus-2}, we get that
\begin{align*}
z_nz_{n+2}-z_{n+1}^{2}
&=
\left(\sum_{k=0}^{n}w_{n,k}\right)
\left(\sum_{k=0}^{n}w_{n,k}\ell_{n,k}\right)
-\left(\sum_{k=0}^{n}w_{n,k}g_{n,k}\right)^2.
\end{align*}
This can be rewritten as
\begin{align*}
z_nz_{n+2}-z_{n+1}^{2}
&=
\left[
\left(\sum_{k=0}^{n}w_{n,k}\right)
\left(\sum_{k=0}^{n}w_{n,k}g_{n,k}^{2}\right)
-\left(\sum_{k=0}^{n}w_{n,k}g_{n,k}\right)^2
\right]\notag\\
&\quad+
\left(\sum_{k=0}^{n}w_{n,k}\right)
\left(\sum_{k=0}^{n}w_{n,k}
\widetilde{\ell_{n,k}}\right),
\end{align*}
where $\widetilde{\ell_{n,k}}=\ell_{n,k}-g_{n,k}^{2}$ for $0\leq k\leq n$.
Since $w_{n,k}\geq0$, by Lemma~\ref{cauchy-schwarz},
the expression in square brackets is nonnegative.
Therefore, it remains to prove that
\begin{align*}
\left(\sum_{k=0}^{n}w_{n,k}\right)
\left(\sum_{k=0}^{n}w_{n,k}
\widetilde{\ell_{n,k}}\right)
\geq0.
\end{align*}
For $0\leq k\leq n$, by the definitions of $g_{n,k}$ and $\widetilde{\ell_{n,k}}$, we find
\begin{align*}
\widetilde{\ell_{n,k}}
=&\alpha_{n,k}g_{n+1,k}+\beta_{n,k}\rho_k g_{n+1,k+1}-g_{n,k}^{\,2}\\
=&\alpha_{n,k}\bigl(\alpha_{n+1,k}-\alpha_{n,k}\bigr)\nonumber+\Bigl[
\alpha_{n,k}\bigl(\beta_{n+1,k}-\beta_{n,k}\bigr)
+\beta_{n,k}\bigl(\alpha_{n+1,k+1}-\alpha_{n,k}\bigr)
\Bigr]\rho_{k}\nonumber\\
&+\beta_{n,k}\bigl(\beta_{n+1,k+1}-\beta_{n,k}\bigr)
\rho_{k}\rho_{k+1}+\beta_{n,k}^{2}\rho_{k}(\rho_{k+1}-\rho_{k}),
\end{align*}
which is nonnegative, since $\alpha_{n,k},\beta_{n,k},\rho_k,\rho_{k+1}$, and
\begin{align*}
\alpha_{n+1,k}-\alpha_{n,k},\,
\beta_{n+1,k}-\beta_{n,k},\,
\alpha_{n+1,k+1}-\alpha_{n,k},\,
\beta_{n+1,k+1}-\beta_{n,k},\,
\rho_{k+1}-\rho_k    
\end{align*}
are nonnegative by the assumptions of the theorem. Thus, $\widetilde{\ell_{n,k}}\geq 0$.
We now consider the case when $x_j=0$ for some $j\geq0$.
For $k\geq1$, since $x_k^2\leq x_{k-1}x_{k+1}$, we have $x_{k}=0$ if $x_{k-1}x_{k+1}=0$. 
Applying this observation repeatedly, we obtain $x_k=0$ for every $k\geq 1$. In particular, if $x_0=0$, then clearly $z_n=0$ for every
$n\geq0$, and hence $(z_n)_{n\geq0}$ is log-convex. It remains to consider the case $x_0>0$ and $x_k=0$ for $k\ge 1$. In this case, $z_n=T_{n,0}x_0$ for every $n\geq0$. For $n\ge 1$, by the recurrence \eqref{T_recurrence}, we have 
\begin{align*}
T_{n-1,0}T_{n+1,0}-T_{n,0}^2&=T_{n-1,0}(\alpha_{n,0}T_{n,0})-(\alpha_{n-1,0}T_{n-1,0})T_{n,0}\\
&=\bigl(\alpha_{n,0}-\alpha_{n-1,0}\bigr)T_{n-1,0}T_{n,0}\ge 0,
\end{align*}
where the last inequality follows from the assumptions of Theorem~\ref{log-convex}.
Thus $(T_{n,0})_{n\geq0}$ is log-convex, and so is $(z_n)_{n\geq0}.$ This completes the proof.
\end{proof}

\section{Strong $q$-Log-Convexity Preservation}\label{Strong $q$-Log-Convexity Preservation}

In this section, we aim to prove
Theorem~\ref{thm:strong-q-log-convexity-preserving}.
Throughout this section, let $T=\{T_{n,k}\}_{n,k\geq0}$ be the
nonnegative triangular array defined by \eqref{T_strong}, with
$(a_n)_{n\geq0}$, $(b_n)_{n\geq0}$ and $(c_n)_{n\geq0}$ satisfying the
assumptions of Theorem~\ref{thm:strong-q-log-convexity-preserving}. As will be shown later, 
the row generating polynomials 
\begin{align*}
P_n(t)=\sum_{k=0}^{n}T_{n,k}t^k, \qquad n\geq 0
\end{align*}
will play an important role in our proof of Theorem~\ref{thm:strong-q-log-convexity-preserving}. 
Before giving the details, we briefly explain the main idea of the proof.
Given a strongly $q$-log-convex sequence
$(X_k(q))_{k\geq0}$ in $\mathbb R_{\geq0}[q]$, let
\begin{align*}
Z_n(q)=\sum_{k=0}^{n}T_{n,k}X_k(q),\qquad n\geq0.
\end{align*}
Our aim is to show that for any $n,m\geq 0$ the difference
\begin{align*}
\mathcal{T}_{n,m}(Z;q)=Z_n(q)Z_{n+m+1}(q)-Z_{n+1}(q)Z_{n+m}(q)
\end{align*}
is a polynomial in $q$ with nonnegative coefficients. 
For this purpose, define a linear map
$\Phi:\mathbb R[u,v]\longrightarrow\mathbb R[q]$ by
\begin{align*}
\Phi(u^iv^j)=X_i(q)X_j(q),\qquad i,j\geq0.
\end{align*}
By setting
\begin{align*}
D_{n,m}(u,v)
=P_n(u)P_{n+m+1}(v)-P_{n+1}(u)P_{n+m}(v),
\end{align*}
we obtain
\begin{align*}
\mathcal{T}_{n,m}(Z;q)=\Phi(D_{n,m}(u,v)).
\end{align*}
Thus, we need to show that  $\Phi\bigl(D_{n,m}(u,v)\bigr)$ has nonnegative coefficients for every
$n,m\geq0$. By definition we will have $\Phi\bigl(f(u,v)\bigr)=0$ if $f(u,v)=-f(v,u)$.
On the other hand, the strong $q$-log-convexity of $(X_k(q))_{k\geq0}$ implies that for any $i,j\geq 0$
both $\Phi(u^iv^j)=X_i(q)X_j(q)$ 
and
\begin{align*}
\Phi((v-u)u^iv^{i+j})=X_i(q)X_{i+j+1}(q)-X_{i+1}(q)X_{i+j}(q)
\end{align*}
have nonnegative coefficients as polynomials in $q$. 
The strong $q$-log-convexity of $(Z_n(q))_{n\geq0}$ will immediately follow if we can express
$D_{n,m}(u,v)$ as a nonnegative linear combination of polynomials $u^iv^j$ and $(v-u)u^iv^{i+j}$ together with some antisymmetric terms. 

As indicated above, we need to find a suitable expansion of $D_{n,m}(u,v)$. To this end, we first reformulate the 
recurrence \eqref{T_strong} as 
\begin{align}\label{eq:row-polynomial-recurrence}
P_{n+1}(t)=(a_n+b_nt)P_n(t)+c_nt(1-t)P_n'(t),
\qquad n\geq0,
\end{align}
with $P_0(t)=1$.
For notational convenience, define the operator
\begin{align}\label{eq-operator-rnt}
R_n^{(t)}:=(a_n+b_nt)+c_nt(1-t)\partial_t
\end{align}
and hence \eqref{eq:row-polynomial-recurrence} can be rewritten as $P_{n+1}(t)=R_n^{(t)}(P_{n}(t))$.
One can verify that
\begin{align}\label{eq:Dnm-operator}
D_{n,m}(u,v)
=\bigl(R_{n+m}^{(v)}-R_n^{(u)}\bigr)
  \bigl(P_n(u)P_{n+m}(v)\bigr).
\end{align}
This inspires us to study the expansion of $P_n(u)P_{n+m}(v)$.
We first consider the case of $m=0$. 
Recall that Wang and Yeh~\cite{Wang-Yeh-2005} obtained the following result. 

\begin{lem}[{\cite[p.~69]{Wang-Yeh-2005}}]\label{lem:real-rootedness} 
Let $f(t)$ be a nonzero polynomial with nonnegative
real coefficients and only real zeros, and let
$\alpha,\beta,\gamma,\delta\in\mathbb R$.
If $\alpha\delta\geq\beta\gamma$, then
\begin{align*}
(\alpha t+\beta)f(t)+t(\gamma t+\delta)f'(t)
\end{align*}
is either zero or has only real zeros.
\end{lem}

Based on lemma \ref{lem:real-rootedness} and the recurrence \eqref{eq:row-polynomial-recurrence}, we have the following factorization of $P_n(u)P_n(v)$.

\begin{lem}\label{lem:row-polynomial-factorization}
For any $n\geq0$ such that $P_n(t)\not\equiv0$, let
$d_n=\deg P_n(t)$.
Then the product $P_n(u)P_n(v)$ can be factorized as
\begin{align*}
P_n(u)P_n(v)
=\kappa_n^2
\prod_{i=1}^{d_n}
\bigl(p_{n,i}^2+p_{n,i}(u+v)+uv\bigr),
\end{align*}
\end{lem}

\begin{proof}
Let us first show that $P_n(t)$ is either zero or has only real
zeros for $n\geq0$. We use induction on $n$. Note that
$P_0(t)=1$. Suppose that the assertion holds for $P_n(t)$.
If $P_n(t)\equiv0$, then the recurrence
\eqref{eq:row-polynomial-recurrence} gives
$P_{n+1}(t)\equiv0$. Otherwise, $P_n(t)$ is a nonzero real-rooted polynomial.
By \eqref{eq-thm3-condition}, we have
$a_n\geq a_0\geq0$ and $c_n\geq0$. For any integers
$s>r\geq0$, we also have
\begin{align}
b_s-b_r
&=\sum_{k=r}^{s-1}(b_{k+1}-b_k)
 \geq\sum_{k=r}^{s-1}
 \bigl((k+1)c_{k+1}-kc_k\bigr)
 =sc_s-rc_r.
\label{eq-bsr}
\end{align}
In particular, taking $s=n$ and $r=0$ in \eqref{eq-bsr}
and using $b_0\geq0$, we obtain
$b_n\geq b_0+nc_n\geq0$
for $n\geq1$. Hence $a_n$, $b_n$ and $c_n$ are all
nonnegative for every $n\geq0$. Take
\begin{align*}
\alpha=b_n,\qquad
\beta=a_n,\qquad
\gamma=-c_n,\qquad
\delta=c_n.
\end{align*}
Clearly,
$\alpha\delta-\beta\gamma
=c_n(a_n+b_n)\geq0.$
It follows immediately from
Lemma~\ref{lem:real-rootedness} that $P_{n+1}(t)$ is either
zero or has only real zeros. This completes the induction.

Moreover, since $P_n(t)$ has nonnegative coefficients, every
nonzero $P_n(t)$ must be a polynomial with nonpositive real
zeros, say
\begin{align*}
P_n(t)=\kappa_n
\prod_{i=1}^{d_n}(t+p_{n,i}),
\end{align*}
for some $\kappa_n>0$ and $p_{n,i}\geq0$.
Multiplying the factorizations of $P_n(u)$ and $P_n(v)$, we
obtain the desired result.
\end{proof}
The product $P_n(u)P_n(v)$ is a nonnegative linear combination
of the polynomials $(uv)^i(u+v)^j$ for nonnegative integers $i,j$;
this is trivial when $P_n(t)\equiv0$, and otherwise follows from
Lemma~\ref{lem:row-polynomial-factorization}.
Note that
\begin{align*}
P_{n+m}(v)=R_{n+m-1}^{(v)}\cdots R_{n+1}^{(v)}R_{n}^{(v)}(P_{n}(v)).
\end{align*}
By successively applying the operators $R_{n}^{(v)}, R_{n+1}^{(v)}, \ldots, R_{n+m-1}^{(v)}$ to each term $(uv)^i(u+v)^j$ will produce a summation of polynomials of the form $u^iv^{i+d}(u+v)^h$, where  
$i,d,h$ are nonnegative integers.
For convenience of notation we will write $u^iv^{i+d}(u+v)^h$ as $B_{i,d,h}(u,v)$ (or denoted simply by $B_{i,d,h}$). Thus $P_n(u)P_{n+m}(v)$ can be written as a linear combination of $B_{i,d,h}$.
As will be seen in the proof below, such a
term arises only when $h=0$, in which case its coefficient
vanishes.
Let $\mathcal C_{n,m}$ be the set of all finite nonnegative linear combinations of the polynomials
$B_{i,d,h}$ with $i,d,h\geq0$, $i+h\leq n$, and $d\leq m$. We find that the expansion of $P_n(u)P_{n+m}(v)$ satisfies the following stronger property. In what follows, to avoid treating the case $h=0$ separately,
we adopt the convention that any term involving $B_{i,d,-1}$
is to be omitted.

\begin{lem}\label{C_m,n}
For any $n,m\geq0$, we have $P_n(u)P_{n+m}(v)\in\mathcal C_{n,m}.$
\end{lem}

\begin{proof} The case $P_n(t)\equiv0$ is trivial. We may therefore assume that $P_n(t)\not\equiv0$ and proceed by induction on $m$. As mentioned above,
the case $m=0$ follows immediately from Lemma~\ref{lem:row-polynomial-factorization}.
Now suppose that $P_n(u)P_{n+m}(v)\in\mathcal C_{n,m}$. We proceed to show that $P_n(u)P_{n+m+1}(v)\in\mathcal C_{n,m+1}$.
By \eqref{eq:row-polynomial-recurrence}, we have 
$P_n(u)P_{n+m+1}(v)
=R_{n+m}^{(v)}\bigl(P_n(u)P_{n+m}(v)\bigr).$
By induction, it suffices to show that 
\begin{align}
R_{n+m}^{(v)}B_{i,d,h}=((a_{n+m}+b_{n+m}v)+c_{n+m}v(1-v)\partial_v)B_{i,d,h}\in\mathcal C_{n,m+1}\label{eq-rnpmbidh}
\end{align}
whenever $i,d,h\geq0$, $i+h\leq n$, and $d\leq m$.
Recall that $B_{i,d,h}(u,v)=u^iv^{i+d}(u+v)^h$. 
A direct calculation gives
\begin{align*}
\partial_v B_{i,d,h}=(i+d)u^iv^{i+d-1}(u+v)^h+hu^iv^{i+d}(u+v)^{h-1},
\end{align*}
and hence
\begin{align*}
v(1-v)\partial_v B_{i,d,h}&=(i+d)u^iv^{i+d}(u+v)^h+hu^iv^{i+d+1}(u+v)^{h-1}\\
&-(i+d)u^iv^{i+d+1}(u+v)^h-hu^iv^{i+d+2}(u+v)^{h-1},
\end{align*}
which can be further written as
\begin{align*}
v(1-v)\partial_v B_{i,d,h}&=(i+d)B_{i,d,h}+hB_{i,d+1,h-1}
-(i+d)B_{i,d+1,h}-hB_{i,d+2,h-1}.
\end{align*}
Note that
\begin{align*}
B_{i,d+2,h-1}&=v^2B_{i,d,h-1}=(v(u+v)-uv)\cdot (u^iv^{i+d}(u+v)^{h-1})\\
&=u^iv^{i+d+1}(u+v)^{h}-u^{i+1}v^{i+d+1}(u+v)^{h-1}\\
&=B_{i,d+1,h}-B_{i+1,d,h-1}.
\end{align*}
Thus, we get
\begin{align}
v(1-v)\partial_v B_{i,d,h}=&(i+d)B_{i,d,h}-(i+d+h)B_{i,d+1,h}\nonumber\\
&+hB_{i,d+1,h-1}+hB_{i+1,d,h-1}.\label{R3}
\end{align}
Substituting \eqref{R3} into \eqref{eq-rnpmbidh} leads to 
\begin{align*}
R_{n+m}^{(v)}B_{i,d,h}
={}&\bigl(a_{n+m}+c_{n+m}(i+d)\bigr)B_{i,d,h}
+\bigl(b_{n+m}-c_{n+m}(i+d+h)\bigr)B_{i,d+1,h}\\
&+c_{n+m}hB_{i+1,d,h-1}+c_{n+m}hB_{i,d+1,h-1}.
\end{align*}
It remains to show that 
\begin{align*}
a_{n+m}+c_{n+m}(i+d)\geq 0,\, c_{n+m}h\geq 0,\, b_{n+m}-c_{n+m}(i+d+h) \geq 0. 
\end{align*}
The first two inequalities are clear by the assumptions of Theorem \ref{thm:strong-q-log-convexity-preserving}.
To verify the third inequality, take $s=n+m$ and $r=0$ in
\eqref{eq-bsr}. Since $b_0\geq0$, $i+h\leq n$, $d\leq m$,
and $c_{n+m}\geq0$, we obtain
\begin{align*}
b_{n+m}
\geq b_{n+m}-b_0
\geq (n+m)c_{n+m}
\geq (i+d+h)c_{n+m}.
\end{align*}
Moreover, the hypothesis that $B_{i,d,h}\in \mathcal C_{n,m}$ forces that all the polynomials $B_{i,d,h},\,B_{i,d+1,h},\,B_{i+1,d,h-1},\,B_{i,d+1,h-1}$ belong to $\mathcal C_{n,m+1}$.
This implies that 
$R_{n+m}^{(v)}(\mathcal C_{n,m})
\subseteq\mathcal C_{n,m+1}$, and hence
$P_n(u)P_{n+m+1}(v)\in\mathcal C_{n,m+1}$, as desired. 
\end{proof}
Let $\mathcal{K}$ denote the set of all polynomials that can be
written as the sum of an antisymmetric polynomial in $u$ and $v$
and a nonnegative linear combination of polynomials of the forms
$u^iv^j$ and $(v-u)u^iv^{i+j}$. With this notation, our main task is  to show  that
$D_{n,m}(u,v)\in\mathcal{K}.$ 
By the preceding lemma and \eqref{eq:Dnm-operator}, it therefore suffices to show that $(R_{n+m}^{(v)}-R_n^{(u)})B_{i,d,h}\in \mathcal{K}$ for each $B_{i,d,h}\in \mathcal C_{n,m}$. 
This is established in the following lemma.

\begin{lem}\label{lem:formal-order}
For any nonnegative integers $n,m,i,d$ and $h$, we have $(R_{n+m}^{(v)}-R_n^{(u)})B_{i,d,h}
\in\mathcal K$ provided $B_{i,d,h}\in \mathcal C_{n,m}$. 
\end{lem}

\begin{proof}
We claim that $B_{i,d,h}\in\mathcal K$ and $(v-u)B_{i,d,h}\in\mathcal K$ for any $i,d,h\geq 0$. 
It is obvious that for any $i,d,h\geq 0$ both $B_{i,d,h}$ and $(v-u)B_{i,d,0}$ belong to $\mathcal K$. 
To prove the claim, we may assume that $d\geq 1$ since $(v-u)B_{i,0,h}=(v-u)(uv)^i(u+v)^h$, which is clearly antisymmetric in $u$ and $v$ and hence belongs to $\mathcal K$. We proceed by induction on $h$. For $d\geq1$, it can be checked that
\begin{align*}
B_{i,d,h}
=B_{i,d+1,h-1}+B_{i+1,d-1,h-1}.
\end{align*}
By the induction hypothesis, both $(v-u)B_{i,d+1,h-1}$ and $(v-u)B_{i+1,d-1,h-1}$ lie in $\mathcal{K}$, and hence so does $(v-u)B_{i,d,h}$. This completes the proof of the claim. 

From now on we assume that $B_{i,d,h}\in \mathcal C_{n,m}$, which implies that $i+h\leq n$ and $d\leq m$.
To continue our proof, we first use \eqref{eq-operator-rnt} to rewrite $R_{n+m}^{(v)}-R_n^{(u)}$ 
as
\begin{align*}
R_{n+m}^{(v)}-R_n^{(u)}&=(a_{n+m}-a_n)+(b_{n+m}-b_n)v+b_n(v-u)\\
&+(c_{n+m}-c_n)v(1-v)\partial_v+c_n\bigl(v(1-v)\partial_v-u(1-u)\partial_u\bigr).
\end{align*}
To expand $(R_{n+m}^{(v)}-R_n^{(u)})B_{i,d,h}$, we need to study the effect of each of the above five summands on the right hand side on 
$B_{i,d,h}$. For the first three summands, their actions on $B_{i,d,h}$ are immediate. The effect of the fourth summand can be given by \eqref{R3}. For the fifth summand, from
\begin{align*}
v\partial_v B_{i,d,h}=(i+d)B_{i,d,h}+hvB_{i,d,h-1},\quad
u\partial_u B_{i,d,h}
&=iB_{i,d,h}+huB_{i,d,h-1}
\end{align*}
we derive that 
\begin{align}
&\bigl(v(1-v)\partial_v-u(1-u)\partial_u\bigr)B_{i,d,h}\notag\\[5pt]
=&(v\partial_v-u\partial_u)B_{i,d,h}
-(v^2\partial_v-u^2\partial_u)B_{i,d,h}\notag\\[5pt]
=&dB_{i,d,h}+h(v-u)B_{i,d,h-1}
-\bigl((i+d)v-iu\bigr)B_{i,d,h}
-h(v^2-u^2)B_{i,d,h-1}\notag\\[5pt]
=&dB_{i,d,h}+h(v-u)B_{i,d,h-1}
-\bigl(dv+i(v-u)\bigr)B_{i,d,h}
-h(v-u)B_{i,d,h}\notag\\[5pt]
=&dB_{i,d,h}-d\,vB_{i,d,h}
-(i+h)(v-u)B_{i,d,h}
+h(v-u)B_{i,d,h-1}.\label{R2}
\end{align}
Combining \eqref{R3} and \eqref{R2} 
we obtain
\begin{align}\label{eq:operator-difference-B}
(R_{n+m}^{(v)}-R_n^{(u)})&B_{i,d,h}
=\bigl(a_{n+m}-a_n+c_nd
+(c_{n+m}-c_n)(i+d)\bigr)B_{i,d,h}\nonumber\\[5pt]
&+(c_{n+m}-c_n)hB_{i+1,d,h-1}
+(c_{n+m}-c_n)hB_{i,d+1,h-1}\nonumber\\[5pt]
&+c_nh(v-u)B_{i,d,h-1}+\bigl(b_n-c_n(i+h)\bigr)(v-u)B_{i,d,h}\nonumber\\[5pt]
&+\bigl(b_{n+m}-b_n-c_nd
-(c_{n+m}-c_n)(i+d+h)\bigr)vB_{i,d,h}.
\end{align}
It remains to prove the nonnegativity of the involved coefficients in the above expansion. 
For the first four summands, the coefficient nonnegativity follows immediately from \eqref{eq-thm3-condition}.
For the fifth summand and the sixth summand, based on \eqref{eq-bsr} and the conditions that $i+h\leq n$ and $d\leq m$, we obtain 
\begin{align*}
&b_n-c_n(i+h)\geq (b_n-b_0)-c_n(i+h)
 \geq c_n(n-i-h)\geq 0,
\end{align*}
and
\begin{align*}
&b_{n+m}-b_n-c_nd-(c_{n+m}-c_n)(i+d+h)\\[5pt]
&\geq (n+m)c_{n+m}-nc_n-dc_n-(i+d+h)(c_{n+m}-c_n)\\[5pt]
&=(m-d)c_n
 +(n+m-i-d-h)(c_{n+m}-c_n)\geq0,
\end{align*}
as desired.
By the previous claim and \eqref{eq:operator-difference-B}, we have $(R_{n+m}^{(v)}-R_n^{(u)})B_{i,d,h}
\in\mathcal K$. This completes the proof.
\end{proof}

We are now in the position to give a  proof of
Theorem~\ref{thm:strong-q-log-convexity-preserving}.

\begin{proof}[Proof of Theorem~\ref{thm:strong-q-log-convexity-preserving}]
Suppose that $(X_k(q))_{k\geq0}$ is a strongly $q$-log-convex sequence of polynomials with nonnegative coefficients.
By definition, for any $r,s\geq 0$ the difference
$$\mathcal{T}_{r,s}(X;q)=(X_r(q)X_{r+s+1}(q)-X_{r+1}(q)X_{s+r}(q))$$
has only nonnegative coefficients. 
For any $n\geq 0$ let
\begin{align*}
Z_n(q)=\sum_{k=0}^{n}T_{n,k}X_k(q).
\end{align*}
We need to show that $(Z_n(q))_{n\geq0}$ is also strongly $q$-log-convex.
Let $D_{n,m}(u,v)$ and $\Phi$ be defined as in the beginning of this section. 
By \eqref{eq:Dnm-operator}, Lemma~\ref{C_m,n} and
Lemma~\ref{lem:formal-order}, we know that for any $n,m\geq0$ there exist an antisymmetric polynomial $g(u,v)$ and finitely many nonnegative real numbers $a_{i,j}$ and $b_{r,s}$ such that
$$D_{n,m}(u,v)=g(u,v)+\sum_{i,j\geq 0}a_{i,j}u^iv^j+\sum_{r,s\geq 0}b_{r,s} (v-u)u^rv^{r+s}.$$
Applying $\Phi$ to both sides leads to 
\begin{align*}
\mathcal{T}_{n,m}(Z;q)= \sum_{i,j\geq 0}a_{i,j}X_i(q)X_{j}(q)+\sum_{r,s\geq 0}b_{r,s} \mathcal{T}_{r,s}(X;q),
\end{align*}
since $g(u,v)=-g(v,u)$ implies that $\Phi(g(u,v))=0$.
Thus $\mathcal{T}_{n,m}(Z;q)$ has only nonnegative coefficients, as desired. The proof is complete. 
\end{proof}

\section{Applications}\label{sec:applications}

In this section, we apply Theorem~\ref{log-convex}
and Theorem~\ref{thm:strong-q-log-convexity-preserving}
to obtain new results on the preservation of log-convexity
and strong $q$-log-convexity for the $r$-colored Eulerian transformation
and $m$-Stirling Eulerian transformation.
 We  also derive several preservation
results appearing in the literature in a unified manner.
\subsection{$r$-Colored Eulerian transformation}

Fix an integer $r\ge1$. Following Steingr\'{\i}msson~\cite[Definitions~1 and~3]{Ste-1994},
an $r$-colored permutation of $[n]$ is represented by a pair
$(\pi,\mathbf s)$, where
$\pi=\pi_1\cdots\pi_n\in\mathfrak S_n$ and
$\mathbf s=(s_1,\ldots,s_n)\in\{0,1,\ldots,r-1\}^n$.
The number of colored excedances of $(\pi,\mathbf s)$ is defined as
\begin{align*}
\operatorname{exc}_r(\pi,\mathbf{s})
=\#\bigl\{i\in[n]:\pi_i>i
\text{ or }\pi_i=i\text{ and }s_i>0\bigr\}.
\end{align*}
Let $A^{(r)}(n,k)$ denote the number of such permutations having
$k$ colored excedances. Set $A^{(r)}(0,0)=1$ and
$A^{(r)}(n,k)=0$ for $k<0$ or $k>n$. The corresponding
$r$-colored Eulerian polynomials and $r$-colored Eulerian
transformation are defined, respectively, by
\begin{align*}
A_n^{(r)}(q)
   =\sum_{k=0}^n A^{(r)}(n,k)q^k,\qquad
z_n^{(r)}
   =\sum_{k=0}^n A^{(r)}(n,k)x_k,
   \qquad n\geq0.
\end{align*}
The polynomials $A_n^{(r)}(q)$ also arise as specializations of
Brenti's bivariate Eulerian polynomials of type~$B$~\cite{Brenti1994}. Let $\mathfrak{S}_n^B$ denote the set of signed
permutations and define
\begin{align*}
B_n(x,q)
=\sum_{\pi\in\mathfrak{S}_n^B}
x^{\operatorname{des}_B(\pi)}
q^{\operatorname{neg}(\pi)},
\end{align*}
where $\operatorname{des}_B(\pi)$ and $\operatorname{neg}(\pi)$
denote, respectively, the number of type~$B$ descents and the number
of negative entries of $\pi\in\mathfrak{S}_n^B$.  Then
$A_n^{(r)}(q)=B_n(q,r-1)$ (see Ma et al.~\cite[Section~7.3]{Ma-Ma-Yeh-Yeh-2024}).

Steingr\'{\i}msson~\cite[Lemma~16]{Ste-1994} established the following recurrence for $A^{(r)}(n,k)$:
\begin{align*}
A^{(r)}(n+1,k)
&=(rk+1)A^{(r)}(n,k)
 +\bigl(r(n-k+2)-1\bigr)A^{(r)}(n,k-1),
\end{align*}
for $n\ge0$ and $0\le k\le n+1$.
It is easy to verify that, for fixed $r\geq 1$, the triangular array 
$T=\{T_{n,k}\}_{n,k\ge0}$ with $T_{n,k}=A^{(r)}(n,k)$ satisfies the conditions of  
Theorem~\ref{thm:strong-q-log-convexity-preserving}. Precisely, we take
\begin{align*}
a_n=1,\qquad b_n=r(n+1)-1,\qquad c_n=r
\end{align*}
in \eqref{T_strong}. We have the following result. 

\begin{cor}\label{cor:colored-eulerian-preservation}
For every integer $r\ge1$, the $r$-colored Eulerian
transformation preserves strong
$q$-log-convexity and hence log-convexity as well.
\end{cor}

\begin{rem}
We would like to point out that some special cases of Corollary \ref{cor:colored-eulerian-preservation} are also interesting. 
For $r=2$, the $r$-colored Eulerian transformation $z_n=\sum_{k=0}^nA^{(r)}(n,k)x_k$ is called the Eulerian transformation of type $B$ by Liu and Zhu \cite{Liu-Zhu-2015}, who conjectured that if $(x_k)_{k\geq 0}$ is log-convex, so is $(z_n)_{n\geq 0}$. Our result is clearly stronger than Liu and Zhu's conjecture. Another interesting special case of  Corollary \ref{cor:colored-eulerian-preservation} is to take $X_k(q)=q^k$ for $k\geq 0$ in the $r$-colored Eulerian transformation $Z_n(q)=\sum_{k=0}^nA^{(r)}(n,k)X_k(q)$. Since $(q^k)_{k\geq 0}$ is clearly strongly $q$-log-convex, the above corollary implies that $(A_n^{(r)}(q))_{n\geq0}$ is strongly $q$-log-convex. 
It is worth mentioning that the strong $q$-log-convexity of $(A_n^{(r)}(q))_{n\geq0}$ also follows from \cite[Proposition 3.4]{Liu-Zhu-2015} and  \cite[Theorem 15 and Theorem 20]{Ste-1994}.
\end{rem}

\subsection{$m$-Stirling Eulerian transformation}

Fix an integer $m\ge1$.
An $m$-Stirling permutation of order $n$ is a permutation
of the multiset $\{1^m,2^m,\ldots,n^m\}$ such that every
entry between two consecutive occurrences of $i$ is greater
than $i$, for each $i\in[n]$.
The case $m=2$ gives the classical Stirling permutations
introduced by Gessel and Stanley~\cite{Gessel-Stanley-1978}; see also
Janson, Kuba and Panholzer~\cite[Section~2.1]{JKP-2011}
for the general case.

For an $m$-Stirling permutation
$\sigma=\sigma_1\cdots\sigma_{mn}$, set
$\sigma_0=\sigma_{mn+1}=0$ and call $i$ a descent if
$\sigma_i>\sigma_{i+1}$, where $0\le i\le mn$.
In particular, every nonempty permutation has a terminal
descent.
Let $E^{(m)}(n,k)$ count the $m$-Stirling permutations of
order $n$ with $k$ descents.
Set $E^{(m)}(0,0)=1$ and $E^{(m)}(n,k)=0$ for
$k<0$ or $k>n$.
The associated polynomial sequence and transformation are defined,
respectively, by
\begin{align*}
E_n^{(m)}(q)
=\sum_{k=0}^n E^{(m)}(n,k)q^k,\qquad
z_n^{(m)}
=\sum_{k=0}^n E^{(m)}(n,k)x_k,
\qquad n\geq0.
\end{align*}
The standard insertion construction for $m$-Stirling permutations
(see \cite[Section~2.1]{JKP-2011}) yields the recurrence
\begin{align*}
E^{(m)}(n+1,k)
&=kE^{(m)}(n,k)
 +(mn-k+2)E^{(m)}(n,k-1),
\end{align*}
for $n\geq0$ and $0\leq k\leq n+1$.
It follows that, for each fixed $m\geq1$, the triangular array
$T=\{T_{n,k}\}_{n,k\geq0}$ defined by
$T_{n,k}=E^{(m)}(n,k)$ satisfies the conditions of
Theorem~\ref{thm:strong-q-log-convexity-preserving}.
Indeed, in \eqref{T_strong} we may take
\begin{align*}
a_n=0,\qquad b_n=mn+1,\qquad c_n=1.
\end{align*}
Thus we obtain the following corollary.

\begin{cor}\label{cor:m-stirling-log-convex}
For every integer $m\geq1$, the $m$-Stirling Eulerian
transformation preserves strong $q$-log-convexity and hence
log-convexity.
\end{cor}

\begin{rem}
This corollary implies that the sequence
$\bigl(E_n^{(m)}(q)\bigr)_{n\geq0}$ is strongly
$q$-log-convex. In fact, it also has the $q$-SM property,
which follows from a result due to P\'etr\'eolle, Sokal and
Zhu~\cite{Petreolle-Sokal-Zhu-2023}.
\end{rem}

\subsection{Classical triangular transformations}

We conclude by showing how
Theorem~\ref{thm:strong-q-log-convexity-preserving}
recovers two classical preservation results, namely, those for the
binomial transformation and the Stirling transformation of
the first kind. Let $c(n,k)$ denote the  Stirling numbers of
the first kind. They satisfy
\begin{align*}
\binom{n+1}{k}=\binom{n}{k}+\binom{n}{k-1},\quad
c(n+1,k)=n c(n,k)+c(n,k-1),
\end{align*}
where $\binom{0}{0}=c(0,0)=1$, and both $\binom{n}{k}$ and
$c(n,k)$ are taken to be zero whenever $k<0$ or $k>n$. Comparing these recurrences with
\eqref{T_strong}, we obtain
\begin{align*}
(a_n,b_n,c_n)=(1,1,0)
\qquad\text{and}\qquad
(a_n,b_n,c_n)=(n,1,0),
\end{align*}
respectively. In both cases, the hypotheses of
Theorem~\ref{thm:strong-q-log-convexity-preserving} are 
satisfied. The theorem therefore yields the following known
preservation results.
\begin{cor}\label{cor:classical-log-convexity}
The binomial transformation and the Stirling transformation
of the first kind preserve strong $q$-log-convexity and log-convexity.
\end{cor}
This recovers the corresponding results in
\cite{Liu-Li-2016, Liu-Wang-2007,Zhu-Sun-2015}.
Other classical examples can be treated similarly. For instance,
Theorem~\ref{log-convex} also recovers the log-convexity preservation
of the Stirling transformation of the second kind and the unsigned
Lah transformation; see \cite{Liu-2011,Liu-Wang-2007}.

\noindent \textbf{Acknowledgements.} 
This work is supported by the Fundamental Research
Funds for the Central Universities.

\noindent \textbf{AI disclosure.} 
The authors used ChatGPT for preliminary computational exploration, including testing mathematical properties and exploring counterexamples. We reviewed and revised the content as needed and independently verified every step of the mathematical arguments and computations to ensure the correctness of the results. The authors take full responsibility for the content of the published article.

\end{document}